# A Unified View on Planning, Scheduling and Dispatching in Production Systems

Kan Wu

*School of Mechanical and Aerospace Engineering, Nanyang Technological University, Singapore, phone: 65-6790-5584; fax: 65-6792-4062; e-mail: kan626@gmail.com*

# A Unified View on Planning, Scheduling and Dispatching in Production Systems

Planning, scheduling and dispatching play critical roles in the operations of a supply chain. Their definitions are clearly given through a unified view in this paper. The distinction between planning and scheduling is analyzed from the view point of microeconomics and queueing theory. The distinction between scheduling and dispatching is analyzed from the view point of computational complexity and hierarchical decomposition. Based on the elasticity of price and capacity, planning can be separated into demand planning or capacity planning. Scheduling period is the time horizon where price and average production cost are insensitive to the production rate. The critical roles of the master production schedule and move targets in job scheduling have been explained through the concept of hierarchical decomposition. Dispatching is the last layer of job scheduling in the hierarchical decomposition. The advantage of pull and push systems has been compared and analyzed systematically.

## 1   Introduction

The terminologies of production planning, scheduling and dispatching are commonly used in practice and literature. While they are all related to arranging limited resources to meet the production goals, it is well accepted that production planning deals with the resource arrangement in the long run, and scheduling or dispatching focuses on relatively short term or even real time information. However, there is no rigorous distinction among them. Sometimes, those three levels are studied in literature with different names. Bitran and Dasu (1992) classified planning, scheduling and dispatching as strategic, tactical and operational levels and studied their corresponding mathematical models through an empirical and non-integrated manner.

According to Mönch et al. (2013), planning is performed with a time horizon ranging from months to years in the semiconductor industry, scheduling is the process of allocating scarce resources over time, and dispatching is the activity to assign the next job to be processed from a set of waiting jobs. Hopp and Spearman (2008) also declared that the time horizon for long-term planning is around 6 months to 5 years. A few questions arise from those definitions. Why is the planning horizon months to years but not weeks to months? And, what is the distinction between scheduling and dispatching? The conventional approach for scheduling problems is to solve an optimization model with specific objectives, such as makespan, maximum lateness, or total completion time, etc. (Allahverdi et al. 2008, Pinedo 2012). Through the model, if the scheduling problem is solved, there is no need for dispatching, since jobs can be simply dispatched based on the specified sequences. On the other hand, dispatching is commonly a heuristic which considers local and real time information, such as first-in-fist-out, shortest processing time, and earliest due date, etc. (Uzsoy et al. 1994, Wein 1988). Both scheduling and dispatching aim at giving job processing sequences. Hence, if we have one, there is no need for the other. Sometimes, there is no clear distinction between both, and their use can be confusing. For example, Wein (1988) studied the topic of "scheduling semiconductor wafer fabrication" by comparing the performance of various dispatching rules.

If the roles of scheduling and dispatching are nearly identical, can we have a system which only has planning and dispatching without scheduling in practice? In air traffic management, a normal passenger will get at least three departure times for the airplane: The first one is given when they buy the ticket, the second one is given at the check-in counter, and the last one is the actual departure time from the runway. If the time on the ticket belongs to planning and the actual departure time comes from dispatching, what is the role of the second one? If we do not have an updated departure time at the check-in counter, we may arrive at the gate at the wrong time, since the new departure time can be different from the one on the ticket. In the semiconductor manufacturing, in addition to the committed date to the customers and the actual time to process jobs based on dispatching, there is commonly a daily move target to guide the dispatching. Why and when do we need those extra layers between planning and dispatching? What are their relations to scheduling?

It s well known that most scheduling problems are extremely difficult to solve (Garey et al. 1976, Lenstra et al. 1977), so are the scheduling problems in a semiconductor fab (Uzsoy, et al. 1994). But why can a fab perform consistently without severe disruptions most of the time? To eliminate the gap between theory and practice, it is of fundamental importance to investigate the problems from a unified view and define the roles of planning, scheduling and dispatching clearly without ambiguity.

The contributions of this paper are three-fold. First, rigorous definitions of scheduling and planning are given. The distinction between planning and scheduling is analyzed from the view point of microeconomics and queueing theory. Second, the relations between scheduling and dispatching are analyzed from the view point of modeling complexity and hierarchical decompositions. To eliminate the gap between theory and practice, the roles of empirical approaches of master production schedules and move targets in scheduling are introduced. Third, the pull and push concepts in production systems are clarified. The distinct roles of pull and push systems in planning, scheduling and dispatching are explained.

We start from defining basic terminologies and giving major assumptions in Section 2. The distinction between planning and scheduling are analyzed in Section 3. The roles of scheduling and dispatching are given in Section 4. Conclusion is given in Section 5.

## 2 Definition and Assumption

### 2.1 Definition

To serve the purpose of the latter analysis, we have to clearly define the terminologies first, although some of them may seem obvious and commonly used in practice.

According to Pearsall and Hanks (1998) *production* is "the action of making or manufacturing from components or raw materials." To provide a more complete definition, we define production with its objective as follows.

**Definition 1** (Production)

The action of manufacturing products from components or raw materials to meet the demand by limited resources.

Definition 1 not only describes its function, but also specifies the objective. The three key elements in the definition are supply (i.e., components or raw materials), demand, and resources. In other words, production is to match supply with demand by limited resources. The demand here specifically refers to customer requirements of production quantity and sojourn time for a specific product. The resources can be machines or operators in a production system.

Comparing to scheduling and dispatching, planning generally refers to the activities in longer time horizons. According to Perrault (2008), *planning* is "the act or process of making a plan to achieve or do something." To meet the organization goal of production systems, planning is defined as follows.

**Definition 2** (Planning)

The act or process of making a *plan* to optimize the profit of the system.

In terms of production, a (production) *plan* is about how to match supply with demand by limited resources in order to obtain the optimal profit. Although an organization may sometimes sacrifice its profit to achieve some strategic goals, we assume the objective is still to optimize the profit in the long run. According to Perrault (2008), *scheduling* is "to plan something at a certain time." It succeeds the decisions from the planning level, and aims at realizing the decisions at a certain time. From the view point of production, scheduling is defined as follows.

**Definition 3** (Scheduling)

The action of manufacturing products from components or raw materials to realize the planning level decisions by limited resources *at a certain time*.

In order to optimize the profit at a certain time, scheduling has to satisfy the demand (i.e., production quantity and sojourn time) from the planning level. In the optimization models of a scheduling problem, the planned demand is reflected in the objective function, components or raw materials are the work-in-process (WIP), and resources are captured by the constraints. Based on Blackstone et al. (1982), the definition of dispatching is defined as follows.

**Definition 4** (Dispatching)

The activity to assign the next job to be processed from a set of jobs awaiting service.

An essential difference between dispatching and planning is that dispatching makes decisions for *immediate actions*, but planning makes decisions for *future activities*.

Depending on the dispatching policy, a production system can be either classified as a pull or push system. According to Kimura and Terada (1981), in a push system, its production and inventory control is based on the forecast value, while in a pull system, the replenishment of inventory at each stage is ordered by the succeeding process at the rate it has been consumed. According to Karmarkar (1989), "a pull system initiates production as a reaction to present demand, while push initiates production in anticipation of future demand." Along the same lines, we define a pull system as follows.

**Definition 5** (Pull system)

A pull system makes decisions based on present information.

A pull system makes decisions based on the most updated information, such as current shop floor conditions, incoming job arrival times and updated demand information. When making present decisions, it simply uses present information, and the feedback can be used for closed-loop control (Spearman and Zazanis 1992). When making future decisions, it does not use forecast but regards the future as present and uses the most updated information. In contrast with a pull system, a push system is defined as follows.

**Definition 6** (Push system)

Relative to the time instant that a decision is realized, decisions of a push system are made based on past information.

When making present decisions, a push system uses past information. Hence, the control is open-looped. When making future decisions, a push system uses forecast based on the present information. The present information becomes the past information when the decision is realized or implemented.

## 2.2  Major Assumption

Planning, scheduling and dispatching are closely related to each other but have distinct roles in a production system. Their relations will be clarified in Section 3 and 4. In this paper, we specifically study the production systems under the following assumptions:

1. Average production cost decreases in production rate,
2. Price decreases in production rate in the long run,
3. Price is constant in production rate in the short run,
4. Average production cost decreases in mean sojourn time in the long run,
5. Average production cost is constant in mean sojourn time in the short run,
6. Price decreases in mean sojourn time.

The first assumption implies that the product enjoys the economies of scale. Since the production cost consists of both variable and fixed costs, and the fixed cost per job decreases as the production rate increases, the average production cost decreases in production rate.

Under some weak conditions, the second assumption holds in any market in the long run. First, in both monopoly and imperfectly competitive markets, the price decreases as the production rate increases in the long run (Nicholson 1995). Furthermore, based on queueing theory, increasing production rate implies the increase of mean sojourn time under fixed system capacity. Due to the sixth assumption, price decreases in production rate even in a perfectly competitive market. On the other hand, since it takes time to achieve market equilibrium, the price elasticity is infinite in the short run, which justifies the third assumption.

**Remark.** To understand the price-production rate relation in oligopolistic differentiated products markets, Berry et al. (1995) studied the U.S. automobile prices in market equilibrium. In addition to the above justification, price also decreases in production rate when a production system makes multiple products. Since the average selling price depends on the capacity allocation, if the costs of all products are about the same, an optimal capacity allocation will satisfy the demand with a higher price first. Hence, average price decreases as the total production rate increases. Mallik and Harker (2004) investigated the optimal profit through capacity allocations among different products.

The fourth assumption is caused by capacity expansion and can be justified by the queueing theory. Since under the same demand (1) mean sojourn time decreases when system utilization is lower, and (2) lower utilization is achieved by adding more capacity (e.g. purchasing more machines), the mean sojourn time decreases with the higher investment on capacity, thus the higher average production cost. However, since it takes time to ramp up capacity, the average production cost is constant in mean sojourn time in the short run.

**Remark.** It usually takes considerable time to ramp up capacity in practical manufacturing systems. For example, in the semiconductor industry, it is common to take two ~ three years to complete the process of constructing a new facility until it is ready for production. Even installing a new machine to an existing facility may take quarters from purchasing to installation and commissioning.

To gain higher market share and enter the market earlier with less competition, customers usually would like to pay a higher premium for the shorter sojourn time. Hence, if a production system can maintain lower sojourn time than its competitors, it can charge its customers a higher premium in the long run. In the short run, customers commonly have unexpected urgent requests for shorter sojourn times (e.g. hot lots in a semiconductor fab) and are willing to pay a higher premium as well. This price differentiation (due to different sojourn time requirements) is commonly maintained in the contract negotiated by salesmen. Hence, the more expedited jobs the system can adopt without sacrificing the sojourn times of normal jobs, the more profit the system can make. Furthermore, if the sojourn times follow a specific distribution, reducing the mean sojourn time will increase the service level (i.e., the percentage of jobs that meet the committed date), and thus increase customer satisfaction and product value to the customers. Hence, price is decreasing in mean sojourn time in both short and long runs.

Some of the properties regarding long term behaviors are analyzed in the next section in detail.

## 2.3   Long Term Behavior Analysis

Assume that customers always maximize their own utility, and the utility follows the Cobb-Douglas utility function, $U = \Lambda_X^\alpha \Lambda_Y^\beta$, where $\Lambda_X$ is the total production rate of $X$ from all manufacturers, $\Lambda_Y$ is the total production rate of a pseudo product $Y$ which represents the effect of all other products and $\alpha + \beta = 1$ (Cobb and Douglas 1928). Under the *Ceteris paribus* assumption, $\Lambda_X = \alpha V P_Y^\beta P_X^{\alpha-1} / (\alpha^\alpha + \beta^\beta)$, where $P_X$ is the price of $X$, $P_Y$ is the price of $Y$ and $V$ is the optimal utility (Nicholson 1995). Hence, when $P_Y$ is fixed,

$$P_X = (k/\Lambda_X)^{\frac{1}{1-\alpha}}, \tag{1}$$

where $k$ is a constant and $\Lambda_X = \sum_i \lambda_i$. $\lambda_i$ is the production rate of the $i$th manufacturer.

Under a fixed supply $\Lambda_X$ of $X$, customers generally would like to pay a higher price for a shorter lead time. If the lead time is too long, the demand in the market may disappear or be replaced by its substitutes. Hence, the price that a customer would like to pay will drop to zero if the lead time is too long. Eq. (1) should consider the impact from lead time for completeness. Assuming the lead time of $X$ is $l$, based on Eq. (1) and the above discussion, the price curve is modified as follows,

$$P_X = ((a - bl^\gamma)/\Lambda_X)^{\frac{1}{1-\alpha}}, \tag{2}$$

where $a, b, \gamma \geq 0$, and $a - bl^\gamma \geq 0$. The values of $a, b, \gamma$ and $\alpha$ should be determined by empirical data through regression. Based on Eq. (2), when $\Lambda_X$ is fixed, $P_X$ is lower if $l$ is longer as shown in Figure 1 and Figure 2. The price is higher when the total production rate (i.e., supply) is lower.

In general, a product can be either classified as an innovative or a functional product (Fisher 1997). For functional products such as the staples that people buy in retail stores, since the demand is more stable and can be predicted reliably in advance, price would be insensitive to lead time as long as the lead time is not too long. The price decreases slowly in mean lead time at the beginning but decreases more apparently after a threshold as being captured by the curves in Figure 1. On the other hand, for an innovative product such as smart phones or laptops, customers may even wait in a long queue or pay a higher price to have the product earlier. The price would drop rapidly at the early stage as the mean lead time increases. This scenario is captured by the curves in Figure 2. Since the price also drops dramatically when the supply $\Lambda_X$ increases, a considerable portion of the profit of an innovative product is contributed at the early stage of the product life cycle before the emergence of its substitutes from competitors.

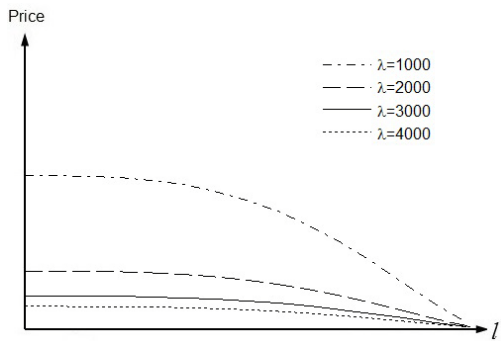

Figure 1 Relations between price and lead time for functional products ($a = 3000, b = 1, \gamma = 3$ and $\alpha = 0.3$)

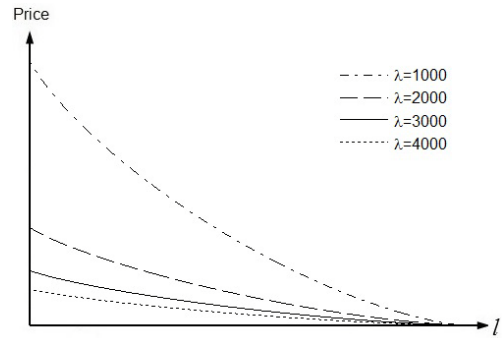

Figure 2 Relations between price and lead time for innovative products ($a = 3000, b = 1, \gamma = 0.8$ and $\alpha = 0.3$)

Eq. (2) reflects the price that a customer would like to pay for each ($\Lambda_X$, $l$) in a market. There is no specification of the relation between $\Lambda_X$ and $l$. Since a supply chain can be viewed as a queueing network and manufacturing sojourn times usually constitute a considerable portion of supply lead time, we use manufacturing sojourn times to represent supply lead time in the following derivation. For a specific manufacturer $i$, the manufacturer will have its own *performance curve*, which characterizes the trade-off between the mean sojourn time and throughput rate of the production system. When all stations consist of a single server, the system mean sojourn time can be approximated as follows (Wu and McGinnis 2012),

$$l_i \cong k_1 \left( \frac{\rho_{BN}}{1 - \rho_{BN}} \right) \frac{1}{\mu_i} + k_2 \left( \frac{\lambda_i}{k_3 - \lambda_i} \right) \frac{1}{k_3} + PT_f, \tag{3}$$

where $l_i$ is mean sojourn time of manufacturer $i$, $PT_f$ is expected total process time, $\mu_i$ is bottleneck service rate (or system capacity), $\lambda_i$ is the production rate of the manufacturer, $\rho_{BN}$ is bottleneck utilization (or $\lambda_i / \mu_i$), $k_1$ is the bottleneck variability, $k_2$ is the variability of a composite station which represents all non-bottleneck stations, and $k_3$ represents the effective capacity of the composite station. For stability, we need to ensure $\lambda_i$ is smaller than $\mu_i$. Note that $\Lambda_X = \sum_i \lambda_i$ in Eq. (2). For systems with multiple server stations, please refer to Wu and McGinnis (2012).

From Eq. (3), we have

$$\lambda_i = \left( -B \pm \sqrt{(B^2 - 4AC)} \right)^+ / 2A, \tag{4}$$

where $A = l_i + \frac{k_1}{\mu_i} + \frac{k_2}{k_3} - PT_f$, $B = PT_f \mu_i + PT_f k_3 - l_i \mu_i - l_i k_3 - \frac{k_1 k_3}{\mu_i} - \frac{\mu_i k_2}{k_3}$ and $C = l_i \mu_i k_3 - PT_f \mu_i k_3$.

Replacing $\lambda_i$ in Eq. (2) by Eq. (4) and assuming $\lambda_j$ in Eq. (2) are fixed for $j \neq i$, the price $P_X$ can be expressed as a function of $l$. Since the curves behave differently in monopoly and competitive markets, they will be discussed separately. If the manufacturer is a monopoly (i.e., $\lambda_i = \Lambda_X$ and $l_i = l$), product price decreases when sojourn time becomes longer. If machine variabilities are all identical, since $l_i$ is a function of $\lambda_i$, $\mu_i$, and $k_i$, the possible curves of $P_X$ vs. $l$ are shown in Figure 3 to Figure 5. For functional (or innovative) products, the price-sojourn time curves are captured by the Price 1 (or Price 2) curves in Figure 3 to Figure 5.

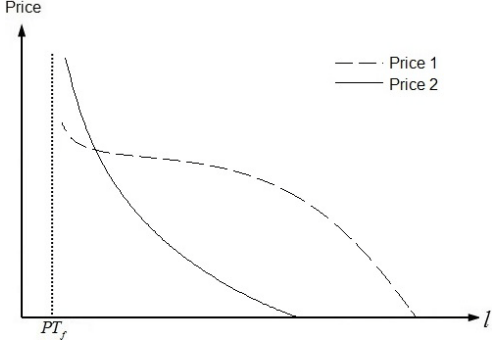

Figure 3 Price curves for a monopoly when $\mu_i$ and $k_i$ are fixed and $\lambda_i$ is changing

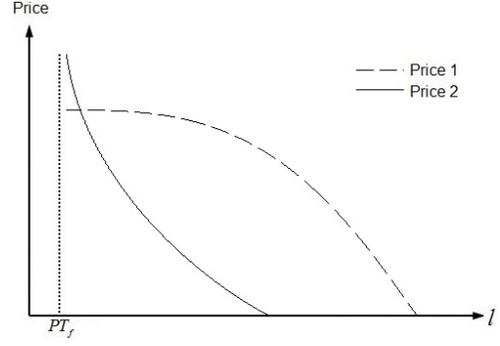

Figure 4 Price curves for a monopoly when $\lambda_i$ and $k_i$ are fixed and $\mu_i$ is changing

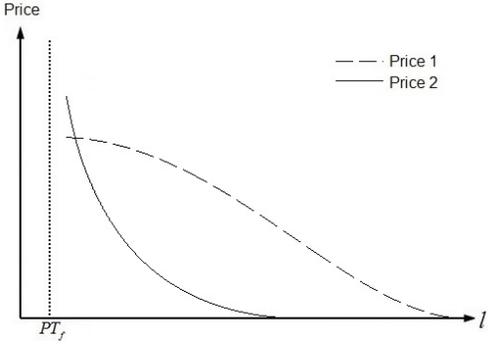

Figure 5 Price curves for a monopoly when $\lambda_i$, $\mu_i$, $k_2$ and $k_3$ are fixed and $k_1$ is changing

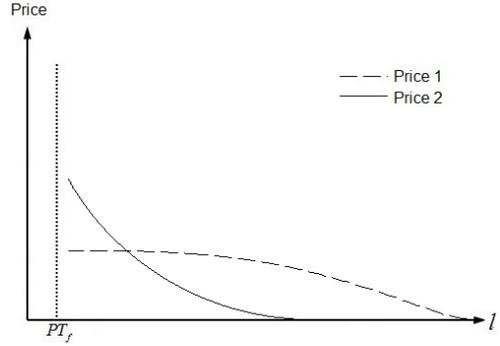

Figure 6 Price curves for a perfectly competitive market

In Figure 3, when $\mu_i$ and $k_i$ are fixed and $\Lambda_X$ (i.e. $\lambda_i$), is changing, with different parameters in Eq. (2), the curve of $P_X$ behaves differently. The Price 1 curve in Figure 3 occurs when there exists $l$ such that $\frac{d^2 P_X}{d\,l^2} = 0$, and the Price 2 curve occurs when $\frac{d^2 P_X}{d\,l^2} > 0$ for all $l$. Note that both Price 1 and 2 curves in Figure 3 diverge to infinity when $l$ approaches $PT_f$. Since under a given capacity $\mu$, based on Eq. (3), the mean sojourn time can be small only when the production rate $\Lambda_X$ (or $\lambda_i$) is small. To achieve economic equilibrium by balancing supply and demand, price has to be very high if $\Lambda_X$ is extremely small as characterized by Eq. (2).

In a perfectly competitive market, no single manufacturer can dominate the whole market (i.e., $\lambda_i$ is relatively smaller than $\Lambda_X$). Under the same settings as the $a, b, \gamma$ and $\alpha$ in Figure 1 and Figure 2, when $\mu_i$ and $k_i$ are fixed and $\lambda_i$ is changing, the price-sojourn time curves of the manufacturer $i$ are shown in Figure 6, where the Price 1 curve refers to the setting in Figure 1 and the Price 2 curve refers to the setting in Figure 2. Since the market is under perfect competition, the production rate of manufacturer $i$ has rare impact on the total production rate $\Lambda_X$. Hence, the curve behaves similar to one of the curves in Figure 1 or Figure 2. Since $\Lambda_X$ is close to a constant, the manufacturer sees the price-sojourn time curve solely depended on customer preference.

When a market is under imperfect competition, the price-sojourn time curve of a specific manufacturer will behave somewhat between the curves of a monopoly and a perfectly competitive market. Cournot, Stackelberg or Bertrand models from game theory could be assumed for the detailed analysis (Nicholson 1995).

Assume that the capacity of station $j$ of manufacturer $i$ is $\mu_{ij}$ jobs per day, unit cost of capacity of station $j$ is $r_{ij}$, variable cost is $v_0$ per job, lifetime of the manufacturing system is $t_i$ and production rate (i.e., production quantity per period) is $\lambda_i$ jobs per day. Hence,

$$\text{Average production cost per job of manufacturer } i = \sum_j r_{ij}\mu_{ij}/\lambda_i t_i + v_0, \tag{5}$$

Replacing $\lambda_i$ in Eq. (5) by Eq. (4), the average production cost can be expressed as a function of $l_i$.

Assuming the product is functional and the manufacturer is a monopoly, based on Eq. (2), (3) and (5), the curves under Assumptions 1 and 2 are depicted in Figure 7 when $\lambda$ is changeable. Based on Eq. (2), (4) and (5), the curves under Assumptions 4 and 6 are depicted in Figure 8 when $\mu$ is changeable.

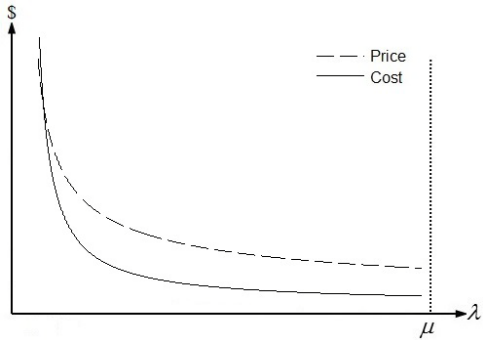
Figure 7 Curves under Assumptions 1 and 2

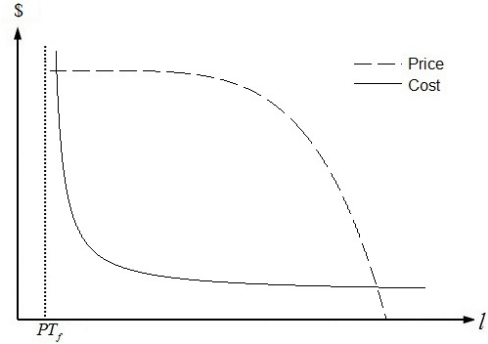
Figure 8 Curves under Assumptions 4 and 6

## 3  Planning and Scheduling

Planning is the process of making a plan to optimize the profit of the system. In the following, we analyze how to optimize the profit of a system under the six assumptions mentioned previously.

The first observation is that due to Assumptions 2 and 3 or Assumptions 4 and 5, system profit can only be fully explored in the long run when selling prices and production costs are variables. Based on the planning time horizon, we could have three types of planning processes: long term planning (or planning), capacity planning and demand planning. Let $t_p$ denote the time horizon that the price is indifferent to the production rate, and $t_c$ denote the time horizon that the production cost is indifferent to the mean sojourn time. If the time horizon for long term planning is $t$, we have the following result.

**Property 1** (Time horizon of long term planning)

$t \geq max(t_p, t_c)$.

When the capacity can be adjusted, we specifically call the planning process capacity planning. If the time horizon for capacity planning is $t$, we have the following result.

**Property 2** (Time horizon of capacity planning)

$t \geq t_c$.

The difference between the two curves in Figure 8 is the profit that a company can make at different mean sojourn times through capacity planning. When the price can be affected by the production rate, we call the planning process demand planning. If the time horizon for demand planning is $t$, we have the following result.

**Property 3** (Time horizon of demand planning)

$t \geq t_p$.

The difference between the two curves in Figure 7 is the profit that a company can make at different production rates through demand planning. However, to obtain the optimal total profit, one has to consider both capacity planning and demand planning at the same time. For the long term planning, based on Eq. (3) and (5), the response surface of production cost, $f_c(\lambda, l)$, is delineated in Figure 9. Based on Eq. (2) and (3), the response surface of price, $f_p(\lambda, l)$, is delineated in Figure 10.

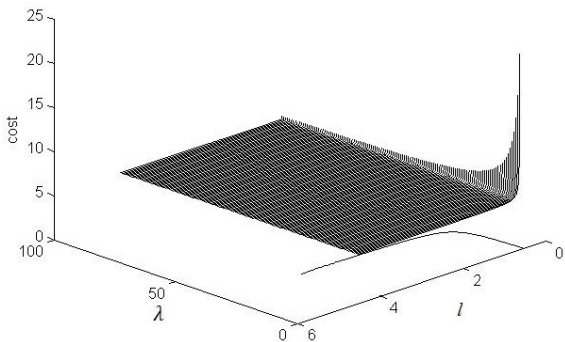

Figure 9 Response surface of production cost

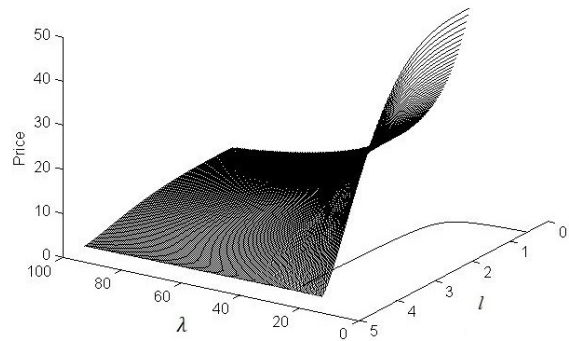

Figure 10 Response surface of price

The two response surfaces in Figure 9 and Figure 10 are combined into Figure 11. Note that not all points on the response surfaces are feasible, since a feasible solution also needs to satisfy Eq. (3), which is the performance curve on the mean sojourn time versus production rate plane. There is a unique performance curve corresponding to each $\mu$.

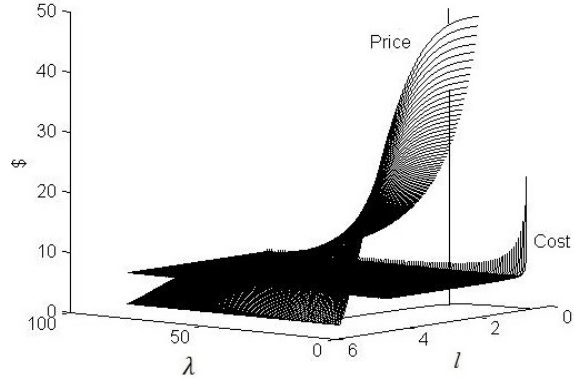

Figure 11 Response surfaces of production cost and price

Let $x$ denote profit, $y$ be the optimal total profit, $p$ be price, $c$ be cost, $\lambda_i$ be production rate and $l_i$ be the mean sojourn time. Since $x = p - c$, the profit of a product is the difference between the two response surfaces along the performance curves in Figure 11. To optimize the total profit, a system should be operated at the point where the profit is maximized, i.e.,

$$\text{Max}_{(\lambda_i, l_i)} \; \lambda_i (p - c)$$
$$s.t. \; p \in f_p(\lambda_i, l_i),$$
$$\quad c \in f_c(\lambda_i, l_i),$$
$$\quad (\lambda_i, l_i) \text{ satisfies Eq. (3).}$$
$$\quad \sum_j r_{ij}\mu_{ij} \leq \text{capital investment of manufacturer } i,$$

where $f_p(\lambda_i, l_i)$ is the response surface of price and $f_c(\lambda_i, l_i)$ is the response surface of production cost. Let $(\bar{\lambda}_l, \bar{l}_l) = \text{augmax}_{(\lambda_i, l_i)} \lambda_i(p - c)$. Hence, $(\bar{\lambda}_l, \bar{l}_l)$ is the optimal operational point of the production system. Note that $l_i$ is a function of $\lambda_i$, $\mu_i$, and $k_i$. If $k_i$ is fixed, $l_i$ can be changed only when either $\lambda_i$ or $\mu_i$ is changed. Hence, finding $(\bar{\lambda}_l, \bar{l}_l)$ is the same as looking for the optimal $(\bar{\lambda}_l, \bar{\mu}_l)$ under fixed $k_i$.

By Definition 3, scheduling is the action of manufacturing products from components or raw materials to realize the planning level decisions by limited resources at a certain time. Since it realizes the planning level decisions (i.e., $(\bar{\lambda}_l, \bar{l}_l)$) within a short time horizon, the price is constant in production rate (Assumption 3) and the average production cost is constant in mean sojourn time (Assumption 5). Hence, capacity and demand plans are fixed at the scheduling level, which is consistent with the assumption made by Bitran and Dasu (1992) for a tactical plan. If the time horizon for scheduling is $t$, we have the following property.

**Property 4** (Time horizon of scheduling)

$0 \leq t \leq min(t_p, t_c)$.

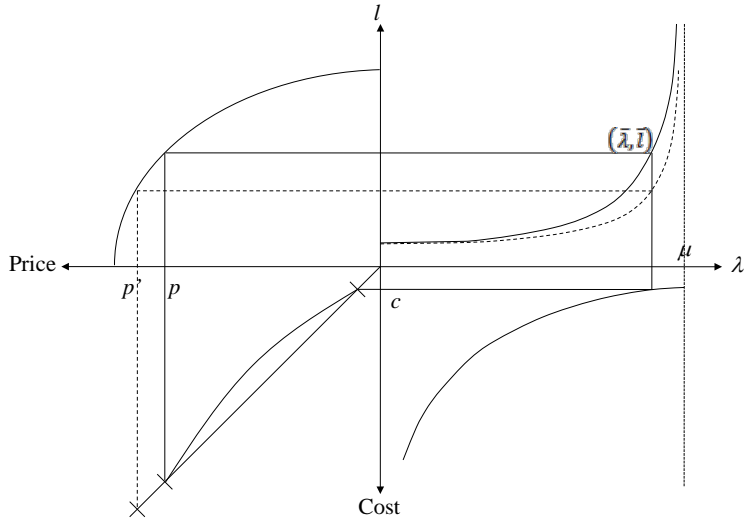

Figure 12 Two-dimensional profit analysis for scheduling

For scheduling, due to Assumption 5 and the planning level decision $(\bar{\lambda}_l, \bar{l}_l)$, the cost-production rate curve is fixed at $l = \bar{l}$ in Figure 9. Due to Assumption 3, the price-sojourn time curve is fixed at $\lambda = \bar{\lambda}$ in Figure 10. Hence, each of Figure 9 and Figure 10 degenerates into a two-dimension graph under Property 4 and we could combine both into one figure as follows. For the clarity of representation, we separate the production cost from the selling price in Figure 11 and project the response surfaces of selling price and production cost on the second and fourth quadrants of a two-dimensional plane in Figure 12.

In Figure 12, the two axes should be read as follows. On the x-axis, the right hand side is the production rate, and the left had side is the selling price. On the y-axis, the upper part is the mean sojourn time and the lower part is the cost per unit product. Since the goal of scheduling is to realize the planning level decisions within a shorter time horizon, note that (1) the optimal operational point $(\bar{\lambda}, \bar{l})$, (2) the price versus mean sojourn time curve (Assumption 6) in the second quadrant, and (3) the production rate versus average production cost curve (Assumption 1) in the fourth quadrant are all succeeded from the planning level decisions.

In Figure 12, the profit of selling a unit of product can be realized in the third quadrant, where $x = p - c$. Although it is not as much degree of freedom as in the planning level, there is still some room to increase the profit at the scheduling level. For example, if one can reduce either the bottleneck or non-bottleneck variability (Wu et al. 2011), based on Eq. (3), the performance curve in the first quadrant will be shifted from the solid line to the dashed line. More profit can be gained by either increasing the production rate or producing more expedited jobs (by reducing the sojourn time of some jobs) while keeping the same production rate. The latter case is expressed by the dash lines in Figure 12.

## 4    Scheduling and Dispatching

Based on Definition 3, scheduling is the action of manufacturing products from components or raw materials to realize the planning level decisions by limited resources within a predefined period satisfying Property 4. To ensure that the planning level decisions on the production quantity and sojourn time can be followed at the scheduling level, the planning level decisions are commonly transformed into a committed date and release date associated with each job at the scheduling level. Assume the quantity of job $i$ committed on the $j$th day is $q_{ij}$.

$$\text{Production quantity} = \sum_j \sum_i q_{ij}, \tag{6}$$

$$\text{Sojourn time} = \text{committed date} - \text{release date} \quad \text{for } \forall \text{ job.} \tag{7}$$

Hence, the planning level decisions on production quantity and sojourn time are equivalent to controlling the committed date and release date of each job. A major difference between these two sets of performance indices is their resolution: the production quantity and sojourn time are aggregated indices and measured in the long run, but the committed date and release date are more detailed and linked to each job. The former two are suitable at the planning level, since they are proposed based on the forecast. The latter two are more detailed and fit the need of shop floor control and scheduling.

Based on Eq. (6) and (7), the objective of scheduling is to manufacture products from raw materials to meet the job committed dates by releasing jobs at the right time and utilizing the limited resources properly over a given time period. Hence, scheduling consists of two major tasks: (i) determining job release dates, and (ii) finding proper allocation between resources and jobs to meet the committed date. Since job release can be viewed as the job-resource allocation at the first step of a process flow, when jobs are released according to the job-resource allocation, (i) is covered by the scope of (ii).

Since the committed date of jobs cannot be satisfied without releasing jobs properly, from the view point of production control, both committed dates and released dates are important. However, from the view point of customers in a supply chain, only the committed date matters. As long as the committed date can be met, customers would not care when the jobs were released. Hence, release dates can be viewed as a means of production control to meet the committed dates. The planning level decisions on the production quantity and sojourn time become to satisfy the committed dates from the view point of customers at the scheduling level.

**Remark.** A key difference between planning and scheduling is that planning deals with demand forecast, while scheduling faces real customer orders. In practice, the demand forecast and customer orders are connected through capacity allocation at the planning level and Available-To-Promise (ATP) at the scheduling level. Capacity allocation allocates capacity to demands according to forecast (Cachon and Lariviere 1999, Mallik and Harker 2004) in the long run and the resolution is usually on the monthly basis (e.g. 1,000 jobs are committed in Jan). Order management at the scheduling level is realized through ATP. The objective of ATP is to provide a response to customer order requests based on resource availability (Zhao et al. 2005). Hence, ATP is the uncommitted portion of inventory or capacity and the resolution is usually on the daily basis.

In production systems, the limited resources can be machines, operators or spare parts. Finding proper allocation between resources and jobs includes the activities such as arranging preventive maintenance schedules, processing bill of material, and assigning flexible machines to machine groups with different capabilities, etc. When the resource specifically refers to a machine group, it is called job-machine scheduling, or simply job scheduling in short (or machine scheduling in some literatures). The allocation between machines and jobs has two directions. In general, allocating jobs to machines occurs more frequently in heavy traffic, while allocating machines to jobs occurs more frequently in light traffic. Job scheduling is realized by controlling job processing sequences at each machine with the consideration of the two-way allocations. An important task of job scheduling is to search for the job processing sequence which guarantees all jobs meet their committed dates. Conventionally, scheduling sometimes simply refers to the job-machine scheduling (Pinedo 2012), but ignores the scheduling activities for other resources. To distinguish both, we call the job-machine scheduling *job scheduling* and the scheduling activities for the other resources *resource scheduling*. Both job scheduling and resource scheduling belong to material requirements planning (or MRP) (Orlicky 1974) and manufacturing resource planning (or MRP II) (Wight 1984). In the following discussion we will mainly focus on job scheduling and clarify its relations with dispatching.

Finding the optimal job processing sequences is difficult in general. Even finding the optimal job processing sequence to achieve the minimum makespan for a simple manufacturing system with only two identical machines (i.e., $P2|C_{max}$) is NP-hard (Garey and Johnson 1975), not to mention more complex problems in practical production systems. Indeed, most scheduling problems in research are NP-hard (Lenstra, et al. 1977). On the other hand, a NP-hard scheduling problem can still be solved within a short period of time if the job size is small. For example, in the clinical trial supply chain, the total setup and process times of an active pharmaceutical ingredient can be one ~ three weeks and there are often less than ten waiting jobs to be scheduled. In this situation, the scheduling problem can be solved efficiently even if it is NP-hard.

Since computation time is the main constraint of job scheduling, it must satisfy Eq. (8) to serve the purpose of scheduling.

$$\text{Computation time} \leq \text{Tolerance time}. \tag{8}$$

The computation time is the duration to get the optimal solution from the mathematical model, and is a function of job size $n$ or $O(f(n))$, i.e., the larger the job size, the longer the computation time. In practical production systems, the job size increases with the number of process steps and scheduling time horizon, i.e., the job size is larger if the process flow or time horizon is longer. Eq. (8) gives a necessary condition to the use of scheduling.

In a stochastic production system, the tolerance time is the duration that the solution remains optimal. Since the solution remains optimal only if the conditions are the same, the tolerance time depends on system variability, which can come from inter-arrival time and service time variabilities (Wu 2005). For example, if the process time of a specific job is longer (or shorter) than predicted, or an unscheduled breakdown occurs, the optimal job sequence may be altered.

In order to ensure true optimality on the sample path basis, job scheduling has to be re-computed whenever the condition changes. Since job shop scheduling problems can be NP-hard, Eq. (8) may not hold especially when the job size is large. On the other hand, since the computation time of a NP-hard problem decreases quickly when the job size is reduced, the decomposition technique is commonly used in practice (Fordyce et al. 1992).

## 4.1 Hierarchical Decomposition

Decomposition is commonly used to deal with complex systems. The purpose is to decompose a complex problem into smaller ones which can be analyzed efficiently. However, to ensure that the decomposed problems achieve the same objective as the original problem, the following two properties have to be satisfied: (1) completeness, and (2) inheritance.

Assume the feasible solution domain of the $i$th subproblem is $D_i$ and the feasible solution domain of the original problem is $D$. Completeness means $D \subset \{D_1 \cup D_2 \cup ...\}$. Hence, the original feasible solution domain is completely covered by the total feasible solution domain of the subproblems. Inheritance means the objective of a subproblem is synchronized with the original problem, so that if the objectives of all subproblems are achieved, the objective of the original problem is also achieved.

Since the computation time of a job scheduling problem is a function of job size, a direct way to satisfy Eq. (8) is to decompose the original problem with a large job size into smaller ones. Because the job size to be scheduled is increasing in the total process steps and scheduling time horizon, to reduce the scheduling job size, one has to either reduce the scheduling time horizon or focus on only a few specific process steps while still satisfying the committed date of each job. Hence, one can either cut the original time horizon into shorter ones (by time) or decompose the process flow into smaller pieces (by space).

Completeness here is satisfied by ensuring all feasible job processing sequences in the original problem also exist in the decomposed subproblems. Since job processing sequences depend on available WIP and processing constraints, completeness can be achieved if no job is missing after decomposition and all constraints between the original problem and subproblems are identical. Hence, when all constraints are succeeded, completeness can be simply satisfied by ensuring that the entire time horizon and process flow are considered and covered by the subproblems. However, the challenge generally comes from converting high level objectives into detailed ones while complying with inheritance.

Since the practical problem size at the scheduling level can be very large and the decomposition process has to satisfy completeness and inheritance, the decomposition of scheduling has to be made insightfully with care. To satisfy Eq. (8) and the two properties, the process is usually realized through multiple layers in practice. Hence, the decomposition process is called hierarchical decomposition.

The first layer of scheduling in the hierarchical decomposition is the master production schedule (MPS). The MPS decomposes monthly capacity allocations into daily required production quantity (i.e., by time) with considering WIP distributions. To manage production lines effectively, MPS has been widely applied to practical manufacturing systems, such as automotive, semiconductor, food packaging and pharmaceutical industries. However, in some manufacturing systems with complex production flows (e.g. semiconductor fabs), the job size can be still too large to satisfy Eq. (8) after the first layer decomposition. Then the second layer decomposition (i.e., by space) has to be imposed. The output of the second layer decomposition is the move target at each process stage, which is used to guide job dispatching and generate job processing sequences.

The purpose of decomposition is to make the job size smaller, so that the optimal job processing sequence can be obtained within a short time frame. When the hierarchical decomposition is imposed, the output from the upper level becomes the objective of its downstream level, and the bottom level gives the job processing sequence for the scheduling problem. Based on Definition 4, dispatching is the activity to assign the next job to be processed from a set of jobs awaiting service. Hence, dispatching is *the last layer of job scheduling in the hierarchical decomposition.* When the job size is already small in the original scheduling problem and no decomposition is needed to obtain the optimal job processing sequence, scheduling is identical to dispatching.

**Property 5** (Dispatching Vs. Scheduling)

Dispatching gives the job processing sequence and is the last layer of job scheduling in the hierarchical decomposition.

### 4.1.1 Master Production Schedule

In practical systems, the scheduling time horizon can be up to months or quarters. Even if all demands have no uncertainty and machines have no unscheduled breakdown, Eq. (8) is unlikely to be satisfied due to the larger number of jobs. In the hierarchical decomposition, the first layer decomposition is realized by decomposing the time horizon through master production scheduling. According to the American Production and Inventory Control Society (Cox et al. 1995), the MPS is "a statement of what the company expects to manufacture ⋯ It represents what the company plans to produce expressed in specific configurations, quantities and dates."

Since the planning level decision is made based on forecast data, it does not reflect current shop floor situations. To generate a feasible plan while satisfying Eq. (8) and complying with inheritance, the planning level decisions have to be decomposed into the objectives for smaller time buckets with the consideration of current shop floor conditions. The MPS translates the planning level decisions of production quantity and sojourn time into the daily required production quantity.

According to demand realization, an MPS can be divided into two portions: the near term portion contains the realized demand (or orders) with committed dates. The remaining is the long term portion and deals with the unrealized demands through ATP. To comply with inheritance, both portions aim at satisfying the committed date from the planning level. The generation of MPS can be expressed as,

(4.1)  $\mathrm{Min} \sum_i H_i$

   *s.t.* capacity and other resource constraints where slackness $H_i$ = (forecast completion date of job $i$ – committed date of job $i$)$^+$. Coinciding with the production goal, the MPS matches supply with demand by limited resources, where demand is the job with committed date, supply is the existing WIP and other raw materials, and resources refer to the limited machine capacity. To satisfy inheritance, the job completion date should be equal or earlier than the committed date. Let $q_{ij}$ be the quantity of job $i$ which will be completed on the $j$th day by forecast.

Through Eq. (6), the daily required production quantity of an MPS can be computed. The objective of dispatching becomes to satisfy the MPS daily production quantity.

Since manufacturing systems are stochastic in general, the model settings are often stochastic as well, and the optimal solution is given in terms of expectation in the long term average sense. Hence, the optimal solution from the stochastic optimization model may not be optimal on each individual sample path. For example, due to the transient behavior among time buckets, one cannot exclude the possibility that the total production quantity could be increased by sacrificing the daily production quantity on a specific day but increasing the production quantities on the following days. Because the exact inheritance can be hardly achieved in a stochastic production environment, the decomposition process inevitably leads to sub-optimality in terms of each individual sample path.

**Property 6** (Inheritance in a stochastic production system)
In a stochastic production system, exact inheritance on the sample path basis can be hardly achieved.

In practical manufacturing systems, to reduce the computational complexity, the completion date is generated by considering the on-hand inventory, existing WIP distribution and system capacity in an empirical manner (Proud 2012). Without considering the detailed job processing sequence and machine status, the job completion date is simply computed based on the existing WIP distribution, historical manufacturing sojourn time and a rough estimate of the system capacity (e.g. 100 jobs per day).

Empirically the daily required production quantity of the MPS can be estimated iteratively from the present to the future as follows,

$$Q_j = \min (A_j, D_j, \mu_j), \tag{9}$$

where $Q_j$ is the required production quantity, $A_j$ is the available WIP, $D_j$ is the demand and $\mu_j$ is the capacity on the $j$th day ($j = 1$ at present). There are $n$ stages in a process flow and their indices are arranged in the ascending order. The duration of the MPS time bucket is 1 day. Let $W_m$ be the initial WIP at the $m$th stage, $L_m$ be the historical sojourn time of the $m$th stage, and $q_{ij}$ be the quantity of job $i$ committed on the $j$th day. Then,

$$A_j = \sum_{m=n-R_j^U}^{n-R_j^L} W_m, \tag{10}$$

$$D_j = \sum_i q_{ij} + \left(D_{j-1} - Q_{j-1}\right)^+, \tag{11}$$

where

$R_j^U = \text{Arc} \max_i \sum_{m=n-i}^{n} L_m \leq j$,

$R_j^L = \text{Arc} \max_i \sum_{m=n-i}^{n} L_m \leq j - 1$,

$L_m$ or $W_m$ is 0 if $m \leq 0$ and $D_j$ or $Q_j$ is 0 if $j \leq 0$.

Eq. (10) computes the available WIP by summing up the present WIP at the stage and those at the upstream stations which could be completed within the specified duration. Since the present demand is the total quantity committed on the date, the demand of later dates is calculated forward from present.

Compared with the monthly capacity allocation, an MPS gives detailed daily required production quantity through consideration of the existing WIP distribution. Based on the daily required production quantity, the time horizon of a job scheduling problem can be decomposed into multiple smaller subproblems with daily time buckets.

As the guidance for production control, the time horizon of the near term MPS is required to cover all unfinished jobs on the shop floors. Depending on how early the order arrives, the near term portion can be months. For example, in the semiconductor industry, due to the long manufacturing sojourn time, customer orders usually come a few months earlier than the committed dates. Hence, the time horizon of the near term MPS can be up to months. In addition to production control, the MPS is also used to plan the usage of raw materials and spare parts, where the long term portion plays a critical role in this aspect. Hence, an MPS is also an important component of MRP.

### 4.1.2    Move Target

In semiconductor fabrication facilities, the process flows of a product can be more than a thousand process steps. When the process flow is very long, due to the larger number of WIP, Eq. (8) may not be satisfied even if the time horizon of the original problem has been decomposed into daily time buckets. In this situation, the second layer of the hierarchical decomposition has to be carried out.

The output of the second layer decomposition is move targets. The computation of move targets translates the daily required production quantity from the MPS into the daily production quantity at each process stage by decomposing the process flow into smaller segments. A process stage consists of a few consecutive process steps of the process flow. Those process steps usually work on the same process function and are geographically close to each other.

Due to the long process flow, move targets have been extensively used as a means of production control in the semiconductor industry. In contrast with their importance in practice, move targets do not capture much attention from researchers in academia. Prior literatures are mainly in the forms of conference proceedings (Ham et al. 2006, Lee et al. 2008, Lee et al. 2007, Wu et al. 1998) or patents (Chen and Tsai 2007, Chuang and Lin 2003). The generation of move targets can be expressed as,

(4.2)  Min $\sum_j G_j$

*s.t.* capacity and other resource constraints at each stage where throughput shortage on the $j$th day $G_j$ = (daily required production quantity on the $j$th day – forecast production quantity at the last stage on the $j$th day)$^+$. The move target is the forecast production quantity at each stage in the model (4.2). Coinciding with the production goal, the move target matches supply with demand by limited resources, where demand is the daily required production quantity from the MPS, resources refer to the machine capacity at each stage, and supply is the available WIP at each stage. The available WIP includes the initial WIP at the stage and those at the upstream stations which will arrive at the stage within the MPS time bucket (i.e., one day). To avoid long computation, the daily required production quantity and available WIP at a stage are computed based on historical sojourn time in an empirical manner (Wu, et al. 1998). To satisfy inheritance, the move target at each stage should help realize the required production quantity from the MPS. However, similar to the MPS, due to the stochastic events among stages, the decomposition process inevitably introduces errors into the solutions and results in sub-optimality.

Without considering the detailed job processing sequence and machine status, the move target at the present day can be estimated iteratively from the last to the first stage as follows,

$$M_k = \min\big(A_1^k, D_1^k, \mu_1^k\big), \tag{12}$$

where $M_k$ is move target, $A_1^k$ is available WIP, $D_1^k$ is demand and $\mu_1^k$ is capacity at the $k$th stage on the first day (or present). There are $n$ stages in a process flow and their indices are arranged in the ascending order. The duration of the MPS time bucket is one day. Let $W_i$ be the initial WIP at the $i$th stage, $L_m$ be the historical sojourn time of the $m$th stage, and $P_j$ be the required production quantity from the MPS on the $j$th day. Then,

$$R_k = \text{Arc}\max_i \textstyle\sum_{m=k-i}^{k} L_m \leq 1, \tag{13}$$

$$A_1^k = \textstyle\sum_{i=k-R_k}^{k} W_i, \tag{14}$$

$$D_1^k = P_j + \big(D_1^{k+1} - M_{k+1}\big)^+, \text{ if } R_j^L < n - k \leq R_j^U, \tag{15}$$

where

$R_j^U = \text{Arc}\max_i \sum_{m=n-i}^{n} L_m \leq j$,

$R_j^L = \text{Arc}\max_i \sum_{m=n-i}^{n} L_m \leq j-1$,

$L_m$ if $m \leq 0$, $W_i$ is 0 if $i \leq 0$ and $D_1^{k+1}$ or $M_{k+1}$ is 0 if $k+1 > n$.

Eq. (14) computes the available WIP by summing up the initial WIP at the $k$th stage and those at the upstream stations which will arrive at the stage within one day. Eq. (15) estimates the demand based on the following assumptions: (1) the demand of the last few process stages is set based on the present MPS target, and (2) the demand of the other stages is no less than the MPS target $P_j$ on a later day according to the historical sojourn times $L_m$. Along with the delinquent demand, the demand of all stages can be calculated backward from the last process stage $n$ iteratively.

## 4.2 Dispatching

After we obtain the daily move target of each stage, the job processing sequence (or dispatching) at the $m$th stage and on the $j$th day can be determined through solving the following mathematical program.

(4.3)  Min $\sum_m G_j^m$

> *s.t.* capacity and other resource constraints at each stage where throughput shortage at the $m$th stage and on the $j$th day $G_j^m$ = (move target at the $m$th stage and on the $j$th day – forecast production quantity at the $m$th stage and on the $j$th day)$^+$. Since the objective complies with inheritance (in the long term average sense), the objective of the original problems could be achieved if the objectives of the subproblems are secured. Furthermore, the job size of each subproblem is much smaller than the original problem. Rather than finding the optimum, the objective of dispatching aims at finding feasible job sequences which can meet the move targets. When there exist multiple production schedules with zero objective values (i.e., multiple feasible plans), one may select the one with the least slackness in order to satisfy customer needs.

Note that the original goal of scheduling is to meet the committed dates of jobs under stochastic shop floor conditions. Due to the problem complexity, move targets are used to serve the purpose of hierarchical decompositions. However, an implicit assumption of the optimization model (4.3) is that all jobs could meet their committed dates if the move targets are satisfied. Since the shop floor conditions are stochastic, a certain percentage of jobs would be slower than expected. Those delayed jobs should be accelerated by giving a higher priority. Other than satisfying move targets, the dispatching objective for those jobs should be to reduce their slackness.

Without hierarchical decompositions, a scheduling problem would be too complicated to be solved efficiently. Hence, heuristics are often used and dispatching was studied in the form of heuristics in literature (Blackstone, et al. 1982, Panwalkar and Iskander 1977, Sarin et al. 2011). In this situation, the dispatching policy at the first stage (i.e., job release) can be different from the others and determined by a separate heuristic called job release policies (Fowler et al. 2002, Glassey and Resende 1988, Ragatz and Mabert 1988). This is consistent with the two major tasks of scheduling as we mentioned at the beginning of Section 4, where determining job release dates refers to the job release policy, and finding proper allocation between resources and jobs refers to the dispatching policy at the other stages except for the first one.

According to Definitions 5 and 6, dispatching can be either pull or push. Most dispatching heuristics are push, since they use past information such as current job arrival times, historical service times or historical sojourn times to determine job processing sequences. Among the vast dispatching heuristics, Kanban (Monden 1983) and CONWIP (Spearman et al. 1990) are two well-known pull heuristics. In a Kanban system, cards at a downstream workstation are used to authorize production of the upstream workstation. Since the most updated WIP information is used, a Kanban system is a pull system. A CONWIP system is a closed system, which authorizes job release according to job departures and keeps the total number of jobs in the system as a constant. Since the most updated job departure information is used to determine the job release, a CONWIP system is also a pull system. However, since both Kanban and CONWIP systems use current information partially and selectively, they are pull-oriented heuristics.

Spearman and Zazanis (1992) showed that WIP is easier to control than throughput. Hence, in order to satisfy the objective of the optimization model (4.3), an ideal pull system should determine the dispatching sequence at a workstation based on the most updated WIP distribution at the entire stages, machine status, incoming job arrival times and production target.

The production goal (i.e., forecast production quantity) in the model (4.3) comes from the higher hierarchy of the scheduling level. Since the production goal from the MPS is determined based on the past information (comparing to the time when the dispatching decision is made), the scheduling levels above dispatching are push systems. However, when the production system is deterministic (i.e., no variability from demand and service time, etc.), because all future information can be derived exactly and becomes a subset of current information, a push system reduces to a pull system. The entire production system from planning, scheduling to dispatching can be pull.

**Property 7** (Push reduces to pull)
In a deterministic production environment, a push system reduces to a pull system, and the entire system can be pull.

Just-in-Time (JIT) (Ohno 1982) is a well known philosophy of pull production. The objective of JIT is to minimize the waste by producing items only when needed (Groenevelt 1993), where Kanban plays a critical role in the implementation of JIT. According to Monden (1983), the conditions for Kanban to work well are: a stable product mix, short setups, proper machine layout, job standardization, improvement and automation. Those conditions aim at changing a stochastic production environment into a deterministic one. Indeed, the objective of JIT can only be realized in a deterministic system, which is consistent with Property 7.

It is nice to have a deterministic production system, since the exact inheritance can be secured and true optimality can be obtained in the deterministic sense on each sample path. However, in the real world, the future events most likely will be stochastic. Regulating the behavior of customers or suppliers is not always a wise choice and may be harmful to an organization. Hence, an efficient production system has to strike a balance between regulating customer or supplier behavior and learning how to deal with uncertainty in this imperfect world. Due to the limit imposed by Eq. (8), hierarchical decomposition is a common approach used in stochastic production systems.

When the production environment is stochastic, the corresponding parameters are random variables and the optimization model (4.3) is stochastic. Let the standard deviation of an estimated parameter be $\sigma$, and the objective value of the optimization model (4.3) is $O$.

**Assumption A1** (Better information, better results)

If $\sigma_i \leq \sigma_j$, $O_i \leq O_j$.

For example, the optimal dispatching policy should consider machine status and incoming job arrival times. If a job indeed arrives at 2:00, but was predicted to arrive at 2:10 from a better forecast and at 3:00 from a poor forecast, Assumption A1 claims that the dispatching polity based on the better forecast will give better objective values (e.g. less throughput shortage) than the one based on the poor forecast.

In a stochastic production environment, planning and scheduling have to be done by forecasts. Assume the actual realization occurs at $t_0$, and forecast is made at $t_1$, where $t_0 > t_1$. Comparing to the realization at $t_0$, let the standard deviations of an estimated parameter in the forecast be $\sigma_1$.

**Assumption A2** (Forecast is not prefect)

$0 = \sigma_0 \leq \sigma_1$.

Assumption A2 states that forecasts cannot be always accurate. When making future decisions, parameters in the future have to be estimated based on the current information. Comparing the estimated parameters with the actual data after realization, (1) let the standard deviation be $\sigma_f$ if the estimated parameter is based on forecast, and (2) let the standard deviation be $\sigma_c$ if the estimated parameter is the same as the present information (without forecast). Assumption A3 assumes the forecast is effective.

**Assumption A3** (Forecast is effective)

$\sigma_f \leq \sigma_c$.

For example, if the price of a specific stock is $10 today and we want to estimate the stock price on a later day, Assumption A3 claims that the price based on the forecast should be more accurate than simply using $10. Based on the above assumptions, we have the following properties.

**Property 8** (Pull at dispatching)
A pull system performs no worse than a push system at dispatching.

Based on Definition 4, dispatching decides the next job to be processed from a set of jobs awaiting service. Based on Assumption A2, a pull system, which makes decision based on present information, has more accurate information. Based on Assumption A1, a pull system performs no worse than a push system.

**Property 9** (Hybrid production systems)
In a stochastic production environment, a hybrid production system which combines both pull and push systems will perform better than a pure pull or push system. Specifically, the hybrid production system is pull at the dispatching level and push at the planning and scheduling levels.

Property 9 is based on the following reasoning: (i) At planning and scheduling levels, based on Assumption A3, a push system, which makes decision based on forecast, is more accurate. Based on Assumption A1, a push system performs no worse than a pull system. (ii) Based on Property 8, at the dispatching level, a pull system performs no worse than a push system. Due to (i) and (ii), a hybrid production system performs better than a pure pull or push system.

In a stochastic production system, a pure pull or push system is not the best choice. On the other hand, although purely deterministic systems can be hardly achieved, variance reduction should always be a goal for an organization. Only in deterministic systems, the exact inheritance can be secured, Eq. (8) can be satisfied and true optimality can be obtained on each sample path. This is also the ultimate goal of pull production systems.

## 5 Conclusion

Planning, scheduling and dispatching play critical roles in supply chain management. In this paper, we give clear definition to each from the view point of microeconomics, queueing theory and modeling complexity. Based on the elasticity of price and capacity, planning can be separated into demand planning or capacity planning. Scheduling period is the time horizon where price and average production cost are insensitive to production quantity. According to the type of resources, scheduling can be distinguished into resource scheduling and job scheduling. Dispatching is the last layer of job scheduling in the hierarchical decomposition and would perform better if it is a pull system.

Conventionally, microeconomics only captures the relation between price and quantity. To clarify the role between planning and scheduling, the original price-quantity model has been extended to incorporate sojourn times. For production cost, relation between quantity and sojourn time has been analyzed through the view point of queueing theory.

To solve a job scheduling problem efficiently, practical problems are commonly decomposed into smaller subproblems through hierarchical decompositions. The first layer decomposition decomposes the problem by time and the second layer further decomposes the problem by space. Namely, the MPS decomposes the original problems into smaller subproblems with a daily time bucket. The move target decomposes the MPS into even smaller subproblems by process stages.

Although scheduling problems are often NP-hard, practical production systems are usually operated efficiently. In solving a scheduling problem, researchers often prove it is NP-hard first and then solve the problem through heuristics. However, a scheduling problem in practice is usually solved through hierarchical decompositions. In addition to proving a scheduling problem is NP-hard, it may be equally, if not more, important to study how to solve a scheduling problem efficiently through insightful decomposition. The concept of hierarchical decomposition presented in this paper is just a start. In addition to taking the empirical approach, more systematic and rigorous study on this topic is expected in the future.